\documentclass[twoside,12pt]{amsart} 
\usepackage[dvips]{graphicx}
\DeclareGraphicsExtensions{ps,eps}
\makeatletter
\def\serieslogo@{}
\makeatother
\makeatletter 
\def\@setcopyright{}
\makeatother

\newtheorem{Theorem}{Theorem}[section]

\newtheorem{Proposition}[Theorem]{Proposition}

\theoremstyle{definition}

\theoremstyle{remark}

\begin{document}

\title
{
Solution by convex minimization
of the Cauchy problem for
hyperbolic systems of conservation laws with convex entropy
}

\author{Yann Brenier}
\address{
CNRS, Centre de math\'ematiques Laurent Schwartz,
Ecole Polytechnique, FR-91128 Palaiseau, France.
}
\curraddr{}
\email{yann.brenier@polytechnique.edu}

\maketitle
\markboth{}{}

\noindent
\section*{Abstract}
We show that, for first-order systems of conservation laws with a strictly convex entropy,
in particular for the very simple so-called "inviscid" Burgers equation,
it is possible to address the Cauchy problem by a suitable convex minimization
problem, quite similar to some problems arising in optimal transport or variational mean-field game theory.
In the general case, we show that smooth, shock-free, solutions can be recovered
on some sufficiently small interval of time. In the special situation of the Burgers equation, we further
show that every "entropy solution" (in the sense of Kruzhkov)
including shocks, can be recovered, for arbitrarily long time intervals.

\subsubsection*{Keywords}
Partial differential evolution equations, calculus of variations, conservation laws,
entropy solutions, convex optimization, optimal transport, mean-field games.

\section*{Introduction}
Solving Cauchy problems by convex minimization techniques is definitely not a new idea, in particular in the framework of
linear evolution PDEs, as illustrated by the classical least square method. For instance, in the
case of a linear transport equation such as $\partial_t u+\partial_x u=0$, one can try to
minimize 
$$
\int\int (\partial_t u+\partial_x u)^2 dxdt,\;\;\;u(0,\cdot)=u_0,
$$
where $u_0$ is the initial condition. This typically leads to (degenerate) space-time elliptic problems.
In the framework of nonlinear equations, similar strategies can be used but usually lead to
non-convex ill-conditioned minimization problems. In this paper, we discuss
a somewhat different strategy for the special class of
"systems of conservation laws with convex entropy" \cite{Da,Ma}, 
namely systems of 
first order evolution PDEs of form
\begin{equation}\label{system}
\partial_t U^\alpha+\partial_i (\mathcal F^{i\alpha}(U))=0,\;\;\alpha=1,\cdot\cdot\cdot,m,
\end{equation}
(with implicit summation on repeated indices),
where $U=U(t,x)\in\mathcal W\subset\mathbb{R}^m$, 
$t\in [0,T]$, 
$x\in\mathcal D=(\mathbb{R}/\mathbb{Z})^d$,
$\partial_t=\frac{\partial}{\partial t}$,
$\partial_i=\frac{\partial}{\partial x_i}$, 
$\mathcal{W}$ is a smooth convex open set, while
the so-called
"flux" $\mathcal F:\mathcal W\rightarrow \mathbb R^{d\times m}$ is a smooth function enjoying the symmetry property
\begin{equation}\label{symmetry}
\forall i\in\{1,\cdot\cdot\cdot,n\},\;\;\forall\beta,\gamma\in\{1,\cdot\cdot\cdot,m\},\;\;
\partial^2_{\alpha\beta}\mathcal E
\partial_{\gamma}\mathcal F^{i \alpha}
=\partial^2_{\alpha\gamma}\mathcal E
\partial_{\beta}\mathcal F^{i \alpha},
\end{equation}
for some smooth function, called "entropy",
$\mathcal E:\mathcal W\rightarrow \mathbb{R}$, which is supposed to be 
convex in the strong sense that 
$(\partial^2_{\alpha\beta}\mathcal E)$ is a positive definite matrix, everywhere on $\mathcal W$.
\\
The symmetry condition (\ref{symmetry}) enforces the conservation of entropy, in the sense
that every smooth solution $U$ to (\ref{system}) must satisfy the extra-conservation law
\begin{equation}\label{entropy}
\partial_t (\mathcal E(U))+\partial_i (\mathcal Q^i(U))=0,
\end{equation}
where the "entropy-flux" function $\mathcal Q:\mathcal W\rightarrow\mathbb R^d$ depends
on $\mathcal F$
and $\mathcal E$.
Indeed, (\ref{symmetry}) is equivalent to
$$
\partial_\gamma(\partial_{\alpha}\mathcal E\partial_{\beta}\mathcal F^{i \alpha})
=\partial_\beta(\partial_{\alpha}\mathcal E\partial_{\gamma}\mathcal F^{i \alpha})
$$
which means that $\partial_{\alpha}\mathcal E\partial_{\beta}\mathcal F^{i \alpha}$
is the gradient of some smooth function $Q^i:\mathcal W\rightarrow \mathbb{R}$, namely
$
\partial_{\alpha}\mathcal E\partial_{\beta}\mathcal F^{i \alpha}=\partial_\beta Q^i.
$
So, every smooth solution $U$ of system (\ref{system}) satisfies
$$
-\partial_t (\mathcal E(U))=\partial_\alpha\mathcal E(U)\partial_i (\mathcal F^{i\alpha}(U))
=\partial_\alpha\mathcal E(U)\partial_\beta\mathcal F^{i\alpha}(U)\partial_i U^\beta
=\partial_\beta \mathcal Q^i(U)\partial_i U^\beta
=\partial_i (\mathcal Q^i(U)),
$$
which exactly is (\ref{entropy}).
\\
This important class of PDEs contains many classical models in
continuum mechanics and material sciences (Euler equations of hydrodynamics, Elastodynamics with polyconvex energy,
ideal Magnetohydrodynamics, etc...\cite{Da,Ho,Ma,Se}). 
The simplest example, is, of course, the celebrated and so-called "inviscid" Burgers equation
\begin{equation}\label{burgers}
\partial_t u+\partial_x(\frac{u^2}{2})=0,\;\;\;u\in\mathbb{R},\;\;\mathcal E(u)=\frac{u^2}{2},
\end{equation}
A richer example is the Euler equation of isothermal fluids:
\begin{equation}\label{Euler}
\partial_t \rho+\nabla\cdot q=0,\;\;\partial_t q+\nabla\cdot(\frac{q\otimes q}{\rho})+\nabla\rho=0,
\;\;\mathcal E
=\frac{|q|^2}{2\rho}+\rho\log\rho,\;\;\;\rho> 0,\;\;q\in \mathbb{R}^d.
\end{equation}
Under mild additional conditions, PDEs of that class are (locally) well-posed [typically in
Sobolev spaces $H^s$ for $s>d/2+1$ \cite{Da,Ma}]. Interestingly enough,
in most cases, smooth solutions develop singularities, called "shock waves", in finite time, and they
cease to be continuous. Thus, it is interesting to consider weak solutions $U$, for which the initial
value problem with initial condition $U_0$ means
\begin{equation}\label{weak}
\int_{[0,T]\times\mathcal D}\partial_t W_\alpha U^\alpha+\partial_i W_\alpha F^{i\alpha}(U)
+\int_{\mathcal D}W_\alpha(0,\cdot)U_0^\alpha=0,
\end{equation}
for all smooth functions $W=W(t,x)\in\mathbb{R}^m$ such that $W(T,\cdot)=0$.
Weak solutions are not expected to satisfy the extra-conservation law (\ref{entropy}). In particular
$$
t\in [0,T]\rightarrow \int_{\mathcal D}\mathcal E(U(t,\cdot))
$$
is not expected to be constant in $t$. Whenever a weak solution satisfies
the "entropy inequality" in the sense of distributions
\begin{equation}\label{entropineq}
\partial_t (\mathcal E(U))+\partial_i (\mathcal Q^i(U))\le 0,
\end{equation}
it is called a "weak entropy solutions". The concept of weak  solutions is quite faulty 
since,  for a fixed initial condition $U_0$,
weak solutions may not be unique and the conservation of entropy is generally not true.
(This is now well established in the case of the Euler equations, in Hydrodynamics, through the results of Scheffer, Shnirelman,
De Lellis-Sz\'ekelyhidi Jr.  \cite{Sc,Sh,DS,Wi}.)
However, a weak solution $\tilde U$ cannot differ from a smooth solution $U$ given on 
$[0,T]\times \mathcal D$ as long as 
$$
\int_{\mathcal D}\mathcal E(\tilde U(t,\cdot))\le \int_{\mathcal D}\mathcal E(U(0,\cdot)),\;\;\;a.e.\;t\in [0,T],
$$
which is certainly true for weak entropy solutions, because of (\ref{entropineq}).

The main idea of our paper is to look, given an initial condition $U_0$,
for weak solutions that minimize over $[0,T]$ the time integral
of their entropy. As explained above, this problem is not void since
weak solutions may not be unique and do not conserve their entropy.
Using the trial functions $W$ in (\ref{weak}) as Lagrange multipliers, we get
the following min-max problem
\begin{equation}\label{infsup}
I=\inf_{U}\sup_{W}
 \int_{[0,T]\times \mathcal D} E(U)-\partial_tW_{\alpha}U^\alpha-\partial_iW_{\alpha}F^{\alpha i}(U)
-\int_\mathcal{D}W_\alpha(0,\cdot)U_0^\alpha,
\end{equation}
where $W=W(t,x)\in\mathbb{R}^m$ are smooth functions, vanishing at $t=T$.
This indeed amounts to looking for a weak solution $U$ of our system of conservation laws with initial condition $U_0$ that minimizes
the time integral of its entropy.
Let us now exchange the infimum and the supremum in the definition of $I$
and get the lower bound
\begin{equation}\label{supinf}
J=\sup_{W}\inf_{U}
 \int_{[0,T]\times \mathcal D} E(U)-\partial_tW_\alpha U^\alpha-\partial_iW_\alpha F^{i\alpha}(U)
-\int_\mathcal{D}W_\alpha(0,\cdot)U^\alpha_0
\end{equation}
which can be reduced to the concave maximization problem
\begin{equation}\label{max}
J=\sup_{W}
 \int_{[0,T]\times \mathcal D} -K(\partial_t W,DW)
-\int_\mathcal{D}W_\alpha(0,\cdot)U^\alpha_0
\end{equation}
where $W$ is still subject to $W(T,\cdot)=0$ 
and $K$ is the $convex$ function defined by
\begin{equation}\label{K}
K(A,B)=\sup_{V\in\mathcal W} A_{\alpha}V^\alpha+B_{\alpha i}F^{i\alpha}(V)-\mathcal E(V)
,\;\;\;A\in\mathbb{R}^m
,\;\;\;B\in\mathbb{R}^{m\times d}.
\end{equation}
This concave maximization problem is somewhat similar to the Monge optimal mass transport problem with quadratic cost in its 
so-called "Benamou-Brenier"
formulation \cite{BB,AGS,Vi}.
[This is particularly true in the case of the "inviscid" Burgers equation,
as will be seen in section \ref{section3}.]
\\
In the first part of the paper (section \ref{section2})
we establish that smooth solutions of general systems of conservation
laws with convex entropy can be recovered at least for short enough intervals of time.
In the special case of the "inviscid" Burgers equation, the shortness condition is sharp in the
sense that it exactly corresponds to the formation of a shock. Thus it is tempting to investigate
the possibility of recovering solutions beyond the formation of shocks.
At first glance, this seems difficult. Indeed, solutions of a concave space-time maximization problem
are expected to enjoy some (limited) elliptic regularity,
which does not seem compatible with shocks.
We leave this question open in the general case. Nevertheless, in the second part of this paper
(sections \ref{section3} and \ref{section4}), we are able to prove, in the very elementary
case of the "inviscid" Burgers equation, that, indeed, solutions with shocks 
(more precisely, "entropy solutions", in the sense of Kruzhkov,
which are known to be unique
for each given initial condition \cite{Da,Ho,Ma,Se}) 
solve the maximization problem in a suitable sense.
This result is obtained by combining convex duality arguments (as in \cite{BB}),  of optimal
transport type and related to mean-field games (for which we refer to 
\cite{AGS,BB,CG,LL,Vi}), and classical properties of 
entropy solutions and Hamilton-Jacobi equations (for which we refer to 
\cite{Ba,Da,CS,Ho,Ma,Se}).
Finally, let us mention our previous work \cite{Br}
where the Euler equations of incompressible fluids \cite{Li} have been investigated with similar
ideas.
\subsection*{Acknowledgments}
This work has been partly supported by the grant MAGA
ANR-16-CE40-0014 (2016-2020) and partly performed in the framework of the
CNRS-INRIA team MOKAPLAN.
The author thanks 
Nassif Ghoussoub 
for exciting discussions in the summer of 2016 about his theory of "ballistic optimal transport problem"
\cite{Gh}, which were very influential for the present work.
\section{Recovery of smooth solutions}\label{section2}

\begin{Theorem}\label{consistency}
Let $U$ be a smooth solution of the entropic system of conservation laws (\ref{system}).
Then, as long as $T>0$ is not too large, more precisely as long as
\begin{equation}\label{criterion}
\partial^2_{\beta\gamma}\mathcal E(V)
+(T-t)\partial_i\left(\partial_\alpha \mathcal E(U(t,x))\right)
\partial^2_{\beta\gamma}\mathcal F^{i \alpha}(V)
\end{equation}
stays a positive definite matrix, for all $t\in [0,T]$, $x\in\mathcal D$
and $V\in\mathcal W$,
then the maximization problem (\ref{max}) admits 
\begin{equation}\label{solmax}
W_\alpha(t,x)=(t-T)\partial_\alpha\mathcal E(U(t,x)),\;\;\alpha\in\{1,\cdot\cdot\cdot,m\},
\;\;\;(t,x)\in [0,T]\times\mathcal D,
\end{equation}
as solution.
\end{Theorem}
Notice that condition (\ref{criterion}) requires, at $t=0$,
\begin{equation}\label{criterion0}
\partial^2_{\beta\gamma}\mathcal E(V)
+T\partial_i\left(\partial_\alpha \mathcal E(U_0(x))\right)
\partial^2_{\beta\gamma}\mathcal F^{i \alpha}(V)
\end{equation}
to be a positive definite matrix, for all $x\in\mathcal D$ and $V\in\mathcal W$, which restricts the choice of $T$
with respect to the initial condition $U_0$. Observe, however, than condition (\ref{criterion}) gets
less restrictive as $t$ approaches $T$ and even allows a blow-up of 
$\partial_i\left(\partial_\alpha \mathcal E(U(t,x))\right)$ of order $(T-t)^{-1}$.
\\
As a matter of fact, in the very special and elementary case of 
the "inviscid" Burgers equation (\ref{burgers}) with initial condition $u_0$,
condition (\ref{criterion}) reads
$$
1+(T-t)\partial_x u(t,x)>0,\;\;\;\forall t\in [0,T], x\in\mathbb{R}/\mathbb{Z}
$$
and turns out to be equivalent to (\ref{criterion0}), namely
$$
1+Tu'_0(x)>0,\;\;\;\forall x\in\mathbb{R}/\mathbb{Z}.
$$
This exactly means that 
$T$ is smaller than 
$$
T^*=\inf_{x\in\mathcal D}\;\frac{1}{\max\{-u'_0(x),0\}}\in ]0,+\infty],
$$
which is exactly the first time when a shock forms.
So, at least in this very elementary case, $all$ smooth solutions can be recovered from
the maximization problem in just one step.

\subsubsection*{Proof of Theorem \ref{consistency}}
Since $U$ is supposed to be a smooth solution of (\ref{system}), we have
$$
\partial_t U^\alpha+\partial_\beta\mathcal F^{i\alpha}(U)\partial_i U^\beta=0.
$$
Thus $W$ defined by (\ref{solmax}) satisfies
$$
\partial_t W_\gamma-\partial_\gamma\mathcal E(U)
=(t-T)\partial^2_{\alpha\gamma}\mathcal E(U)\partial_t U^\alpha
=-(t-T)\partial^2_{\alpha\gamma}\mathcal E(U)\partial_\beta\mathcal F^{i\alpha}(U)\partial_i U^\beta
$$
which is equal, thanks to symmetry property (\ref{symmetry}), to
$$
-(t-T)\partial^2_{\alpha\beta}\mathcal E(U)
\partial_{\gamma}\mathcal F^{i \alpha}(U)
\partial_i U^\beta
=-(t-T)\partial_i(\partial_{\alpha}\mathcal E(U))\partial_{\gamma}\mathcal F^{i \alpha}(U)
=-\partial_i W_\alpha\partial_{\gamma}\mathcal F^{i \alpha}(U).
$$
Thus, we have obtained
$$
\partial_t W_\gamma+\partial_i W_\alpha\partial_{\gamma}\mathcal F^{i \alpha}(U)-
\partial_\gamma\mathcal E(U)=0,
$$
which precisely means that, at each point $(t,x)$, $V=U(t,x)$
satisfies the first order optimality condition 
in the definition of $K(\partial_t W(t,x),DW(t,x))$ through (\ref{K}), namely
$$
K(\partial_t W(t,x),DW(t,x))=
\sup_{V\in\mathcal W} \partial_t W_{\gamma}(t,x)V^\gamma+
\partial_i W_{\alpha}(t,x)\mathcal F^{i \alpha}(V)-\mathcal E(V).
$$
Meanwhile, condition (\ref{criterion}) tells us, by definition (\ref{solmax}) of $W$, that
$$
\partial^2_{\beta\gamma}\mathcal E(V)
-\partial_i W_\alpha(t,x)\partial^2_{\beta\gamma}\mathcal F^{i \alpha}(V)
$$
is a positive definite matrix for all $(t,x,V)$, 
which means that, for each fixed $(t,x)$,
$$
V\in\mathcal W\rightarrow 
\partial_i W_\alpha(t,x)\mathcal F^{i \alpha}(V)-\mathcal E(V)
$$
is a concave function. So the first order optimality condition we have already obtained for
$V=U(t,x)$ is enough to deduce that
$$
K(\partial_t W,DW)=\partial_t W_{\gamma}U^\gamma+
\partial_i W_{\alpha}\mathcal F^{i \alpha}(U)-\mathcal E(U).
$$
Thus, integrating on $[0,T]\times \mathcal D$ and using that $U$ is solution of
(\ref{system}), we get
$$
\int_{[0,T]\times \mathcal D} 
K(\partial_t W,DW)+\mathcal E(U)
=\int_{[0,T]\times \mathcal D} \partial_t W_{\gamma}U^\gamma+
\partial_i W_{\alpha}\mathcal F^{i \alpha}(U)=
$$
$$
=\int_{\mathcal D} W_{\gamma}(T,\cdot)U^\gamma(T,\cdot)-W_{\gamma}(0,\cdot)U^\gamma(0,\cdot)
=\int_{\mathcal D} -W_{\gamma}(0,\cdot)U_0
$$
since $U_0$ is the initial condition and, by definition (\ref{solmax}), $W(T,\cdot)=0$.
By definition (\ref{max}), the optimal value $J$ of the maximization problem is larger than
$$
\int_{[0,T]\times \mathcal D} 
-K(\partial_t W,DW)-\int_{\mathcal D} -W_{\gamma}(0,\cdot)U_0.
$$
Thus, we have obtained
$$
J\ge
\int_{[0,T]\times \mathcal D} \mathcal E(U).
$$
But, by definition (\ref{infsup}), the right-hand side is certainly larger than $I$
which is an upper bound for $J$. (Indeed $infsup\ge sup\;inf$ is always true.)
We conclude that $I=J$ which shows that $W$ must be optimal for the maximization problem (\ref{max}),
and completes the proof.
\subsubsection*{End of proof}
\section{Recovery of entropy solutions for the "inviscid" Burgers equation}\label{section3}
As already mentioned, the simplest example of entropic conservation law is 
the so-called "inviscid" Burgers equation
(\ref{burgers}) on the torus $\mathcal D=\mathbb{T}=\mathbb{R}/\mathbb{Z}$, namely
$$
\partial_t u+\partial_x(\frac{u^2}{2})=0,\;\;\;u=u(t,x)\in\mathbb{R},
\;\;x\in\mathbb T,
\;\;t\in [0,T].
$$
Let us recall some of its properties.
This equation is Galilean invariant. [Indeed, for each constant $c$, $u(t,x-ct)+c$
is a solution whenever $u$ is itself a solution.]
Thus, it is not a restriction to assume
that $u_0$ has zero mean on $\mathbb T$. 
This allows us to introduce a unique periodic anti-derivative $\phi_0$
with zero mean on $\mathbb T$, so that $u_0=\phi_0'$.
Given $u_0=\phi_0'$, 
this equation admits a unique so-called entropy solution 
\cite{Da,Ho,Ma,Se} explicitly given by
the Hopf-Cole formula
\begin{equation}
\label{hopf-cole}
u(t,x)=\partial_x \phi(t,x),\;\;\;
\phi(t,x)=\inf_{a\in\mathbb R}\;\;\phi_0(a)+\frac{|a-x|^2}{2t}
\end{equation}
$$
=\inf_{a\in [0,1],\;k\in\mathbb Z}\;\;\phi_0(a)+\frac{|a+k-x|^2}{2t}.
$$
(since $\phi_0$ is $\mathbb Z-$periodic).
The entropy solution stays smooth as long as
$$
a\in\mathbb{R}\rightarrow \phi_0(a)+\frac{|a-x|^2}{2t}
$$
admits a single non-degenerate minimum for each $(t,x)$, 
which is equivalent to
$t\phi^"_0+1>0$
and ceases to be true as soon as $t>T^*$
where
\begin{equation}
\label{T*}
T^*=\inf_{a\in\mathbb T}\;\frac{1}{\max\{-\phi^"_0(a),0\}}\in ]0,+\infty].
\end{equation}
Beyond $T^*$, the entropy solution admits discontinuities in $x$, which are called "shocks".
A striking property of entropy solutions, on top of their forward
uniqueness, is their lack of backward uniqueness after shocks form.
Indeed, from formula (\ref{hopf-cole}), we easily deduce that, for each $T>T^*$,
there is another entropy solution $u^T$ which is shock-free before $t=T$
and coincide with $u$ at time $t=T$.
More precisely:

\begin{Proposition}
\label{substitute}
Let $u$ be the unique entropy solution with initial condition 
$u_0={\phi'}_0$. For any fixed $T>0$, let us introduce
${\phi}_0^T$ such that
$a\in\mathbb R\rightarrow 2T\phi_0^T(a)+a^2$
is the convex hull of
$a\in\mathbb R\rightarrow 2T\phi_0(a)+a^2$.
\\
Let $u^T$ be the unique
entropy solution with initial condition $u_0^T={\phi'}_0^T$.
Then, $u^T(t,\cdot)$ is shock-free for $0\le t<T$ and coincides
with $u$ at time $t=T$.
We call $u^T$ the "shock-free substitute" of $u$ on $[0,T]$.
\end{Proposition}
This concept will be crucially used later on.
\\
\\
Let us now move back to 
the maximization problem \ref{max}.
In the case of the Burgers equation,  with initial condition $u_0$, we get
\begin{equation}
\label{J-burgers0}
J=\sup 
\{
\int_{[0,T]\times\mathbb{T}}
-\frac{\partial_t W^2}{2(1-\partial_x W)}
 -\int_\mathbb{T}W(0,\cdot)u_0
 ,\;\;\;W(T,\cdot)=0\}
\end{equation}
where $W=W(t,x)\in\mathbb{R}$ is a smooth function subject to $\partial_x W\le 1$ (and
$\partial_t W=0$ whenever $\partial_x W=1$).
Indeed, according to (\ref{K}), the value of $K$ in (\ref{max}) is given by
\begin{equation}
\label{K-burgers}
K(A,B)=\sup_{v\in\mathbb R} Av+B \frac{v^2}{2}-\frac{v^2}{2}=\frac{A^2}{2(1-B)},
\end{equation}
whenever $B>1$, $K=0$ for $(A,B)=(0,1)$ and $K=+\infty$ otherwise.
Since
$$
 -\int_\mathbb{T}W(0,\cdot)u_0= \int_\mathbb{T}(W(T,\cdot)-W(0,\cdot))u_0
 =\int_{[0,T]\times\mathbb{T}}\partial_t W u_0
$$
(because $W(T,\cdot)=0$ and $u_0$ does not depend on $t$), the maximization problem
equivalently reads
\begin{equation}
\label{J-burgers}
J=\sup
\{\int_{[0,T]\times\mathbb{T}}
-\frac{\partial_t W^2}{2(1-\partial_x W)}+\partial_t W u_0,\;\;\;W(T,\cdot)=0\}.
\end{equation}
Our goal is to prove that we may recover all entropy solutions of the "inviscid" Burgers equation
from this maximization problem,
in the sense:
\begin{Theorem}\label{shock}
Let $u_0=\phi_0'$ where $\phi_0$ is a smooth function on $\mathbb T$ and both
$u_0$ and $\phi_0$ have zero mean.
Let $u$ be the unique entropy solution of the "inviscid" Burgers equation (\ref{burgers}) with 
initial condition $u_0$ and let $u^T$ be its  "shock-free subsitute"
$u^T$ on $[0,T]$ (as in Proposition \ref{substitute}).
\\
Then, for $any$ fixed time $T>0$, the maximization problem 
(\ref{J-burgers})
admits a generalized solution $W$
for which $q=\partial_t W$ and $\rho=1-\partial_x W\ge 0$ are bounded Borel measures
on $[0,T]\times\mathbb T$, where $\rho(T,\cdot)=1$ and $q$ is absolutely continuous with respect to $\rho$, with density $v=u^T$. 
In particular $v(T,\cdot)=u(T,\cdot)$.
\end{Theorem}

\subsubsection*{Remark}

Notice that $v=u^T$ is a shock-free entropy solution on $[0,T[$.
This is not so surprising, since $v$ is obtained from a concave maximization problem in both space
and time, which is, in some vague sense, a (degenerate) space-time elliptic problem so that
some partial regularity of its solution should be expected.

\subsubsection*{Remark: an analogy with... mountain climbing!}
Through Theorem \ref{shock}, we have a rather unusual way of solving the Cauchy problem, 
somewhat reminiscent of 
some well-known techniques in mountaineering.
Indeed, we may interpret pursuing the entropy solution beyond $T^*$ after
shocks have formed
as walking along a sharp crest. Through the concave maximization problem, 
we prefer accessing  to the point of the crest (namely $u(T,\cdot)$) 
by climbing from a different initial point (namely $u^T(0,\cdot)$) following a less dangerous
way up to destination (through the shock-free substitute $u^T$)! Of course, this strategy is rather cumbersome if our
goal is to explore each point of the crest, but actually safer (*).
\\
\\
{\footnotesize
(*) For instance, so far,
nobody has succeeded
in crossing, without supplemental oxygen, 
the long ridge linking mount Everest
to mount Lhotse 
at about $8000$ meters above sea-level! 
Many thanks to Thomas Gallou\"et for this information.
}
\section{Proof of Theorem \ref{shock}}\label{section4}
The first step of the proof consists in establishing a suitable generalized framework to
solve the maximization problem (\ref{J-burgers}), relying on generalized solutions $W$
of bounded variation on $[0,T]\times\mathbb T$.
Next, we use standard convex analysis to find the corresponding dual (or rather pre-dual)
minimization problem, which is very similar to the most standard optimal transportation problem,
as in \cite{BB}, and actually corresponds to the most elementary
first-order mean-field game \`a la
Lasry-Lions \cite{LL,CG} (**). This problem can be solved almost explicitly, thanks to the properties
of the "inviscid" Burgers equations such as Proposition \ref{substitute},
which completes the proof.
\\
\\
{\footnotesize
(**) At a computational level, using the augmented Lagrangian method of \cite{BB}, 
these two problems differ just by...two lines of code!}

\subsection{A priori estimates and generalized solutions for problem \ref{J-burgers}}

Let us get few estimates. Trivially, we get $J\ge 0$ (just take $W=0$). Next, for each $r>0$
$$
|\int_{[0,T]\times\mathbb T}\partial_t W u_0|
\le 
\frac{r\sup{u_0}^2}{2}
\int_{[0,T]\times\mathbb T}(1-\partial_x W)
+\int_{[0,T]\times\mathbb T}\frac{\partial_t W^2}{2r(1-\partial_x W)}
$$
(by Cauchy-Schwarz-Young's inequality, using that $1-\partial_x W\ge 0$)
$$
=\frac{rT\sup{u_0}^2}{2}
+\int_{[0,T]\times\mathbb T}\frac{\partial_t W^2}{2r(1-\partial_x W)}
\;\;\;{\rm{(since}}\;\;\int_{\mathbb T}\partial_x W=0).
$$
Thus
\begin{equation}
\label{estimate}
|\int_{[0,T]\times\mathbb T}\partial_t W u_0|\le
\frac{rT\sup{u_0}^2}{2}
+\int_{[0,T]\times\mathbb T}\frac{\partial_t W^2}{2r(1-\partial_x W)}
\end{equation}
which already shows (taking $r=1$) that
$$
J\le  \frac{T\sup{u_0}^2}{2}.
$$
So, we can consider a maximizing sequence $W_n$, $n\ge 1$, such that
$$
 \int_{[0,T]\times\mathbb{T}}
\frac{\partial_t W_n^2}{2(1-\partial_x W_n)}-\partial_t W_n u_0\le -J+n^{-1}.
$$
From estimate (\ref{estimate}), taking $r=2$,
we deduce
$$
 \int_{[0,T]\times\mathbb{T}}
\frac{\partial_t W_n^2}{4(1-\partial_x W_n)}\le -J+n^{-1}+T\sup{u_0}^2\le T\sup{u_0}^2+n^{-1}
$$
(since $J\ge 0$). Thus (again by Cauchy-Schwarz' inequality)
$$
\left(\int_{[0,T]\times\mathbb T}|\partial_t W_n|\right)^2
\le 
\int_{[0,T]\times\mathbb T}(1-\partial_x W_n)
\int_{[0,T]\times\mathbb T}\frac{\partial_t W_n^2}{1-\partial_x W_n}
$$
$$
=T\int_{[0,T]\times\mathbb T}\frac{\partial_t W_n^2}{1-\partial_x W_n}\le 4T^2\sup{u_0}^2+4Tn^{-1}.
$$
Since $\partial_x W_n\le 1$, $\int_\mathbb{T}\partial_x W_n=0$ and $W_n(T,\cdot)=0$, this is already enough to deduce that
$W_n$, up to a subsequence, strongly converges in $L^1([0,T]\times\mathbb T)$ to some
limit $W$ which is necessarily a function of bounded-variation on $[0,T]\times\mathbb T$
(which exactly means that $\partial_t W_n$, $\partial_x W_n$ are bounded Borel measures).
Since $K$, defined by (\ref{K-burgers}), is convex,
we further have, by upper semi-continuity, that
\begin{equation}
\label{Jtilde-burgers}
\tilde J=\int_{[0,T]\times\mathbb{T}}
-\frac{\partial_t W^2}{2(1-\partial_x W)}+\partial_t W u_0\ge J,
\end{equation}
but, a priori, $\tilde J>J$ is possible.
As a function of bounded variation, $W(t,\cdot)$ has a limit at $t\uparrow T$, namely
$W(T-,\cdot)$, in the $L^1(\mathbb T)$ sense. However, it is not clear, a priori, that 
this limit is zero, in spite of the fact that $W_n(T,\cdot)=0$. We need a more precise estimate
to prove it. 
For each smooth function $\psi(x)$ and $(s,t)$ such that
$0\le s\le t\le T$, we have
$$
\int_{\mathbb T}\left(W_n(t,\cdot)-W_n(s,\cdot)\right)\psi
=\int_{[s,t]\times\mathbb T}\partial_tW_n\psi
\le 
\sup |\psi|
\sqrt{ 
\int_{[0,T]\times\mathbb T}\frac{\partial_t W_n^2}{1-\partial_x W_n}
\int_{[s,t]\times\mathbb T}1-\partial_x W_n
}
$$
$$
=\sqrt{|t-s|}
\sup |\psi|
\sqrt{ 
\int_{[0,T]\times\mathbb T}\frac{\partial_t W_n^2}{1-\partial_x W_n}
}
\;\;\;{\rm{(since}}\;\;\int_{\mathbb T}\partial_x W_n=0)
$$
$$
\le 2\sqrt{|t-s|}\sup|\psi|\sqrt{T\sup{u_0}^2+n^{-1}}
$$
(as already established). All these estimates show that the $W_n=W_n(t,x)$ are actually uniformly 
bounded in the space $\mathcal W$ of all bounded-variation functions that are uniformly H\"older continuous (of exponent $1/2$) functions of $t$ valued in the space of distributions on 
$\mathbb{T}$, which is enough to guarantee that they converge  to $W$ in $\mathcal W$
(in the weak-* sense, since the closed bounded convex subset of $\mathcal W$ are weak-*
compact). Hence, $W(T,\cdot)=0$ just follows from $W_n(T,\cdot)=0$.
Thus, it is natural to consider the maximization problem (\ref{J-burgers}) in the generalized class
of bounded-variation functions $W(t,x)$ that are uniformly
continuous in $t$ with values in distributions in $x$
such that $W(T,\cdot)=0$. In this large class, the supremum is certainly attained (as can be seen by
using the same a priori estimates as above).
However, the supremum may a priori exceed the value $J$. Let us now show that this is not the case.

\begin{Proposition}
\label{W-BV}
In definition (\ref{J-burgers0}) (equivalently (\ref{J-burgers})),
the optimal value $J$ is achieved in the generalized class $\mathcal W$
of all bounded-variation functions $W=W(t,x)$
defined
on $[0,T]\times\mathbb{T}$ which are uniformly continuous in $t$ with values in distributions in $x$
and satisfy $\partial_x W\le 1$ and $W(T,\cdot)=0$.
\end{Proposition}
\subsubsection*{Proof of Proposition \ref{W-BV}}

It is enough to show that, given a maximizer $W$ in the large class $\mathcal W$,
there is a corresponding
smooth function  $W_\epsilon$ such that $\partial_x W_\epsilon<1$, $W_\epsilon(T,\cdot)=0$
and
\begin{equation}
\label{result}
\lim\inf_{\epsilon\downarrow 0}
\int_{[0,T]\times\mathbb{T}}
\frac{\partial_t W_\epsilon^2}{2(1-\partial_x W_\epsilon)}-\partial_t W_\epsilon u_0\le 
\int_{[0,T]\times\mathbb{T}}
\frac{\partial_t W^2}{2(1-\partial_x W)}-\partial_t W u_0.
\end{equation}
Let us introduce two nonnegative mollifiers on $\mathbb R$, $\zeta$ and $\gamma$
with respective supports in $[0,1]$
and $[-1,1]$. Then we define
$$
W_\epsilon(t,x)=\int_0^1\left(\int_{-1}^1
W(t+\epsilon\tau,x+\epsilon\xi)\gamma(\xi)d\xi\right) \zeta(\tau)d\tau,
$$
where $W(t,x)$ has been extended by $0$ beyond $t=T$.
By construction $W_\epsilon(T,\cdot)=0$, $\partial_xW_\epsilon<1$, and, by Jensen's inequality,
$$
\int_{[0,T]\times\mathbb{T}}
\frac{\partial_t W_\epsilon^2}{2(1-\partial_x W_\epsilon)}\le 
\int_{[0,T]\times\mathbb{T}}
\frac{\partial_t W^2}{2(1-\partial_x W)},
$$
since $W_\epsilon$ has been defined by local average of $W$.
Furthermore,
$$
\int_{[0,T]\times\mathbb{T}}
\partial_t W_\epsilon u_0
=-\int_{\mathbb{T}}
W_\epsilon(0,\cdot) u_0
=-\int_{\mathbb{T}}\int_0^1\int_{-1}^1u_0(x)
W(\epsilon\tau,x+\epsilon\xi)\zeta(\tau)\gamma(\xi)d\xi d\tau dx
$$
has the same limit (because $W(t,x)$ is uniformly continuous in $t$ with values in 
distributions in $x$) as
$$
-\int_{\mathbb{T}}\int_{-1}^1u_0(x)
W(0,x+\epsilon\xi )\gamma(\xi)d\xi  dx
=-\int_{\mathbb{T}}\int_{-1}^1u_0(x-\epsilon\xi)
W(0,x)\gamma(\xi)d\xi  dx
$$
which, itself, converges to
$$
-\int_{\mathbb{T}}u_0(x)
W(0,x)dx,
$$
since $u_0$ has been supposed to be smooth.
This is enough to deduce (\ref{result}) and conclude the proof of Proposition \ref{W-BV}.

\subsection{Primal and dual formulations of Problem \ref{J-burgers}}
In this next step, we get primal and dual formulations for Problem \ref{J-burgers},
after introducing 
$$\rho=1-\partial_x W\ge 0,\;\;\;q=\partial_t W\in\mathbb{R}.$$
We see that, from definition (\ref{J-burgers}) $J$ can be now written
\begin{equation}
\label{MFG-dual}
J=\sup\{
\int_{[0,T]\times\mathbb{T}}
\;\;-\frac{q^2}{2\rho}+q u_0,\;\;\;(\rho,q)\;s.t.\;\;\rho\ge 0,\;\;
\partial_t \rho+\partial_x q=0,\;\;\rho(T,\cdot)=1
\}
\end{equation}
where the supremum is performed over all bounded Borel measures $(\rho,q)$
over 
$[0,T]\times\mathbb{T}$.
This is very close to an optimal transport problem with quadratic cost, in its "Benamou-Brenier" formulation \cite{BB}. 
At this point, with the help of the Fenchel-Rockafellar duality theorem \cite{Brz},
following \cite{BB},
we may assert

\begin{Proposition}\label{propo-dual}
In the maximization problem (\ref{J-burgers}), rewritten as (\ref{MFG-dual}),
the supremum is achieved by some $(\rho\ge 0,q\in\mathbb R)$ enjoying the following
properties: i) $\rho$ is a non-negative measure, ii) $q$ is absolutely continuous with respect to
$\rho$ with a square-integrable density $v=v(t,x)$ so that
$$
J=\int_{[0,T]\times\mathbb{T}}\;\;-\frac{q^2}{2\rho}+q u_0
=\int_{[0,T]\times\mathbb{T}}\;\;-\frac{1}{2}\rho v^2+\rho v u_0.
$$
In addition, we have the duality result
\begin{equation}
\label{MFG}
J=\inf\{
\int_{\mathbb{T}} -\theta(T,\cdot);
\;\;\theta\;{\rm{s.t.}}\;\;
\partial_t \theta+\frac{1}{2}(\partial_x \theta)^2\le 0,\;
\theta(0,\cdot)=\phi_0
\}
\end{equation}
where the supremum is performed over all smooth functions $\theta$ and $\phi_0$
is the unique function with zero mean on $\mathbb{T}$ such that $\phi'_0=u_0$.
\end{Proposition}

The proof of this rather standard result (see also \cite{CCN}
for similar proofs) is postponed to Appendix 1.
\\
\\
The mixed problem (\ref{MFG-dual},\ref{MFG}) can be interpreted as a $relaxed$ $variational$ formulation of a "mean-field game" \`a la Lasry-Lions \cite{LL,CG},
namely
$$
\partial_t\rho+\partial_x (\rho\partial_x\theta)=0,\;\;\;\partial_t \theta+\frac{1}{2}(\partial_x\theta)^2=0,
\;\;\;\rho(T,\cdot)=1,\;\;\;\theta(0,\cdot)=\phi_0,
$$
where the Hamilton-Jacobi equation $\partial_t \theta+\frac{1}{2}(\partial_x\theta)^2=0$
is relaxed as an inequality. As a matter of fact, solving (\ref{MFG}) is very easy, even simpler than 
solving a standard optimal transport problem! Indeed, from standard Hamilton-Jacobi theory
(cf. \cite{Ba} and \cite{CS} as a recent reference)
$$
\partial_t \theta+\frac{1}{2}(\partial_x \theta)^2\le 0,\;
\theta(0,\cdot)=\phi_0
$$
implies
$$
\theta(t,x)\le\;\;\inf_{a\in\mathbb R}\;\;\phi_0(x)+\frac{|a-x|^2}{2t},\;\;\forall t\in ]0,T]
$$
(or, equivalently,
$$
\theta(t,x)\le\;\;\inf_{a\in [0,1],\;k\in\mathbb Z}\;\;\phi_0(a)+\frac{|a+k-x|^2}{2t},\;\;\forall t\in ]0,T],
$$
since $\phi_0$ is $\mathbb Z-$periodic). 
Thus, by saturating this inequality, we immediately obtain
the optimal value for (\ref{MFG}), namely
$$
J=-\int_{\mathbb T} \inf_{a\in \mathbb R}\{\phi_0(a)+\frac{|a-x|^2}{2T}\}dx,
$$
i.e.
\begin{equation}
\label{optiJ}
J=-\int_{\mathbb T} \phi(T,x)dx,
\;\;\;\phi(t,x)= \inf_{a\in \mathbb R}\;\;\phi_0(a)+\frac{|a-x|^2}{2t},
\end{equation}
where we recognize in $\phi(t,x)$ the anti-derivative of the unique entropy solution with initial condition $u_0=\phi'_0$,
i.e. $\partial_x\phi(t,x)=u(t,x)$. As explained in Proposition \ref{substitute}, we also have
$$
\phi(T,x)=\phi^T(T,x),\;\;\;\phi^T(t,x)=\inf_{a\in \mathbb R}\;\;\phi^T_0(a)+\frac{|a-x|^2}{2t},
$$
where $a^2+T\phi^T_0(a)$ is the convex hull of $a^2+2T\phi_0(a).$
Denoting by $\Omega$ the (closed) set of all points $a\in [0,1]$ where $a^2+2T\phi_0(a)$
coincides with its convex hull 
(and, therefore, $1+T\phi"_0(a)\ge 0$),
we observe that, for each fixed $t\in [0,T]$,
$$
a\in\Omega\rightarrow x=a+t\phi_0'(a)
$$
defines a one-to-one change of variable (with range going from $\Omega$
at time $t=0$ to the full torus at time $T$), under which
$$
a+k={\rm{Argmin}}_{\tilde a\in \mathbb R}\;\;\phi_0(\tilde a)+\frac{|\tilde a-x|^2}{2t},
\;\;\phi^T(t,x)=\phi_0(a)+\frac{|a+k-x|^2}{2t},
$$
for a suitable integer $k\in\mathbb Z$ (that may depend on $(t,x)$). In addition,
$$
\partial_x\phi^T(t,x)=\phi'_0(a)=-\frac{a+k-x}{t},
$$
so that
\begin{equation}
\label{a-x}
\phi^T(t,x)=\phi_0(a)+\frac{t}{2}\phi_0'(a)^2,\;\;\;u^T(t,x)=\phi_0'(a).
\end{equation}
Let us now introduce, for each $t\in [0,T]$,
the probability measure $\rho(t,\cdot)$ defined on $\mathbb T$ by
\begin{equation}
\label{def-rho}
\forall f\in C^0(\mathbb T),\;\;\;
\int_\mathbb T f\;\rho(t,\cdot)=\int_\Omega f(a+t\phi_0'(a))\rho_0(a)da,
,\;\;\;\rho_0(a)=1+T\phi_0^"(a),
\end{equation}
which is just the Lebesgue measure at $t=T$ since
$$
\int_\mathbb T\; f\rho(T,\cdot)
=\int_\Omega f(a+T\phi_0'(a))(1+T\phi_0^"(a))da=\int_\mathbb T f(x)dx.
$$
Thus, we have obtained
$$
\int_\mathbb T \phi^T(T,x)dx=\int_\Omega (\phi_0(a)+\frac{T}{2}\phi_0'(a)^2)\rho_0(a)da
=\int_\Omega (\phi_0(a)+\frac{T}{2}\phi_0'(a)^2)(1+T\phi_0^"(a))da,
$$
which provides an explicit value for
the optimum $J$ of Problem \ref{MFG}, namely
\begin{equation}
\label{valueJ}
J=-\int_\Omega (\phi_0(a)+\frac{T}{2}\phi_0'(a)^2)(1+T\phi_0^"(a))da.
\end{equation}
In the last step of the proof of Theorem \ref{shock}, we introduce, for each $t\in [0,T]$,
a companion measure
$q(t,\cdot)$ for $\rho(t,\cdot)$ (already defined by (\ref{def-rho})):
\begin{equation}
\label{def-q}
\forall f\in C^0(\mathbb T),\;\;\;
\int_{\mathbb T}f\;q(t,\cdot)=\int_\Omega
\phi_0'(a)f(a+t\phi_0'(a))\rho_0(a)da,\;\;\;\rho_0(a)=1+T\phi_0^"(a)
\end{equation}
and want to show that $(\rho,q)$ is optimal for the dual problem (\ref{MFG-dual}).
\\
First, 
we immediately check (from (\ref{def-rho},\ref{def-q}))
$$
\partial_t\rho+\partial_x q=0,\;\;\rho(T,\cdot)=1,
$$
which already show that $(\rho,q)$ is an admissible solution for problem (\ref{MFG-dual}).
Next, we see that $q$ is absolutely continuous with respect to $\rho(t,\cdot)$
and observe that its Radon-Nikodym derivative $v(t,\cdot)$
is nothing but $u^T(t,\cdot)$, the "shock-free substitute"
for $u$ on $[0,T]$.
Indeed, 
through the change of variable $$a\in\Omega\rightarrow x=a+t{\phi'_0}(a),$$
we have already used,
$\rho(t,\cdot)$ is the image of $\rho_0(a)da$ (restricted to $\Omega$)
while $v(t,x)$ is 
just equal to $\phi'_0(a)$ (through its very definition (\ref{def-q})), which is also $u^T(t,x)$,
as seen above in (\ref{a-x}).
\\
Then, we compute
$$
\int_{[0,T]\times\mathbb{T}}\;\;-\frac{q^2}{2\rho}+q u_0
=\int_{[0,T]\times\mathbb{T}}\;\;(-\frac{v^2}{2}+v u_0)\rho
=\int_{[0,T]\times\mathbb{T}}\;\;(-\frac{v^2}{2}+v \phi'_0)\rho
$$
$$
=\int_\Omega
\int_0^T
\left(-\frac{1}{2}\phi_0'(a)^2+\phi_0'(a)\phi_0'(a+t\phi_0'(a))\right)dt\rho_0(a)da
$$
$$
=\int_\Omega
\left(-\frac{T}{2}\phi_0'(a)^2+\phi_0(a+T\phi_0'(a))-\phi_0(a)\right)
\rho_0(a)da
$$
$$
=\int_\Omega
\left(-
\frac{T}{2}\phi_0'(a)^2
-\phi_0(a)\right)
\rho_0(a)da.
$$
$$
{\rm{(Indeed}}\;\;
\int_\Omega
\phi_0(a+T\phi_0'(a))\rho_0(a)da=
\int_\Omega
\phi_0(a+T\phi_0'(a))(1+T\phi_0^"(a))da
=\int_\mathbb{T}\phi_0(x)dx=0.)
$$
This is exactly the value (\ref{valueJ}) just obtained for the optimal value $J$ of both
(\ref{MFG}) and (\ref{MFG-dual}).
We conclude that $(\rho,q)$ defined by (\ref{def-rho},\ref{def-q}) 
is indeed optimal for
(\ref{MFG-dual}) while $q$ is absolutely continuous with respect to $\rho$ with density $v=u^T$.
This completes the proof of Theorem \ref{shock}.

\section{Appendix 1: Proof of Proposition \ref{propo-dual}}

It is convenient to introduce (more or less as in \cite{BB}, see also
\cite{CCN})
the Banach space 
$$E=\{(a,b)\in C^0([0,T]\times\mathbb{T})\times C^0([0,T]\times\mathbb{T})\}$$
equipped with the sup-norm
and the convex functions valued in $]-\infty,+\infty]$
respectively defined by
$$
\Phi(a,b)=0,\;{\rm{whenever}}\;a+\frac{b^2}{2}\le 0,\;\;+\infty\;\;{\rm{otherwise}},\;\;\forall (a,b)\in E
$$
$$
\Psi(a,b)=-\int_{\mathbb{T}}\theta(T,\cdot),\;\;{\rm{whenever}}\;
\;\;a=\partial_t\theta,\;\;b=\partial_x\theta,
\;\;\theta(0,\cdot)=\phi_0,
$$
for some smooth function $\theta:[0,T]\times\mathbb{T}\rightarrow \mathbb{R}$,
and $+\infty$ otherwise. [Notice that $\theta$ is then defined without ambiguity.] 
Observe that there is at least one point $(a,b)\in E$ at which $\Phi$ is continuous while
$\Psi$ is finite. [Take, for instance, $a=-1$, $b=0$.] Since both $\Phi$ and $\Psi$ are convex functions,  
the Fenchel-Rockafellar theorem (as stated in \cite{Brz}, chap. I) asserts that
$$
\inf \{\Phi(a,b)+\Psi(a,b),\;(a,b)\in E\}=\max \{-\Phi^*(\rho,q)-\Psi^*(-\rho,-q),\;(\rho,q)\in E'\}
$$
where $E'$ is the dual of $E$, namely the space of all pairs of real valued
bounded Borel measures $(\rho,q)$ on $[0,T]\times\mathbb{T}$, and $\Phi^*$, $\Psi^*$, are
the Legendre transforms of $\Phi$ and $\Psi$. Let us emphasize that notation "max" is used
to express that the supremum is always achieved on the right-hand side (while the infimum
may not be achieved on the left-hand side).
As in \cite{BB}, we first get
$
\Phi^*(\rho,q)=+\infty,
$
unless $\rho\ge 0$, $q$ is absolutely continuous with respect to $\rho$, with a square 
integrable derivative $v$, in which case
$$
\Phi^*(\rho,q)=
\int_{[0,T]\times\mathbb{T}}\frac{1}{2}\rho v^2.
$$
Next, we find
$
\Psi^*(-\rho,-q)=+\infty,
$
unless $\partial_t \rho+\partial_x q=0$, $\rho(T,\cdot)=1$,
in which case
$$
\Psi^*(-\rho,-q)=-\int_{[0,T]\times\mathbb{T}}q u_0.
$$
[Indeed, by definition of  $\Psi$
$$
\Psi^*(-\rho,-q)
=\sup_\theta \int_{[0,T]\times\mathbb{T}}-\partial_t \theta\rho
-\partial_x \theta q +\int_{\mathbb{T}}\theta(T,\cdot)
$$
(where the supremum is performed over all smooth functions $\theta$ s.t. $\theta(0,\cdot)=\phi_0$)
$$
=\sup_{\tilde\theta} \int_{[0,T]\times\mathbb{T}}-\partial_t \tilde\theta\rho
-\partial_x \tilde\theta q -\phi'_0 q+\int_{\mathbb{T}}\tilde\theta(T,\cdot),
$$
(where we have set $\tilde\theta(t,x)=\theta(t,x)-\phi_0(x)$, which vanishes
at $t=0$, and used that $\phi_0$ is a smooth function with zero mean over $\mathbb{T}$). Thus,
whenever $(\rho,q)$ weakly satisfies $\partial_t \rho+\partial_x q=0$, $\rho(T,\cdot)=1$,
$$
\Psi^*(-\rho,-q)=-\int_{[0,T]\times\mathbb{T}}\phi'_0 q=-\int_{[0,T]\times\mathbb{T}}u_0 q,
$$
and $+\infty$ otherwise, as announced.]
\\
So, we have, on one side
$$
\inf \{\Phi(a,b)+\Psi(a,b),\;(a,b)\in E\}=
\inf\{
\int_{\mathbb{T}} -\theta(T,\cdot);
\;\;\theta\;{\rm{s.t.}}\;\;
\partial_t \theta+\frac{1}{2}(\partial_x \theta)^2\le 0,\;
\theta(0,\cdot)=\phi_0
\}
$$
where the infimum is taken over all smooth functions $\theta$, which is
exactly the optimal value of problem (\ref{MFG}), while, on the other side
$$
\max \{-\Phi^*(\rho,q)-\Psi^*(-\rho,-q),\;(\rho,q)\in E'\}
$$
exactly corresponds to (\ref{MFG-dual}).

So we have proven that (\ref{MFG}) and (\ref{MFG-dual})
are dual to each other, provided (\ref{MFG-dual}) is performed over the class of bounded
Borel measures. 
and $\rho_\epsilon(T,\cdot)=1$, which concludes the proof of Theorem \ref{propo-dual}.


\end{document}